\documentclass[12pt]{amsart}
\usepackage{amssymb,amsmath,amsthm}
\usepackage{a4}
\usepackage{enumitem}
\usepackage{color}
\usepackage[colorlinks=true,
linkcolor=webgreen,
filecolor=webgreen,
citecolor=webgreen]{hyperref}
\definecolor{webgreen}{rgb}{0,0,1}
\definecolor{recrown}{rgb}{1,.2,.6}
\usepackage{fullpage}
\usepackage{orcidlink}

\begin{document}
\newtheorem{theorem}{Theorem}
\newtheorem{corollary}[theorem]{Corollary}
\newtheorem{lemma}[theorem]{Lemma}
\newtheorem{proposition}[theorem]{Proposition}
\theoremstyle{definition}
\newtheorem*{examples}{Examples}
\newtheorem{example}{Example}
\newtheorem{conjecture}[theorem]{Conjecture}
\theoremstyle{theorem}
\newtheorem{thmx}{\bf Theorem}
\renewcommand{\thethmx}{\text{\Alph{thmx}}}
\newtheorem{lemmax}{Lemma}
\renewcommand{\thelemmax}{\text{\Alph{lemmax}}}
\theoremstyle{definition}
\newtheorem*{definition}{Definition}
\theoremstyle{remark}
\newtheorem*{remark}{\bf Remark}
\title{\bf On a factorization result of \c{S}tef\u{a}nescu--II}
\markright{}
\author{Sanjeev Kumar$^{\dagger}$ {\large \orcidlink{0000-0001-6882-4733}}}
\address{$~^\dagger$ Department of Mathematics, SGGS College, Sector-26, Chandigarh-160019, India}
\email{sanjeev\_kumar\_19@yahoo.co.in}
\author{Jitender Singh$^{\ddagger,*}$ {\large \orcidlink{0000-0003-3706-8239}}}
\address{$~^\ddagger$ Department of Mathematics, Guru Nanak Dev University, Amritsar-143005, India}
\email{jitender.math@gndu.ac.in}
\date{}
\maketitle
\parindent=0cm
\footnotetext[3]{$^{,*}$Corresponding author: jitender.math@gndu.ac.in; sonumaths@gmail.com}
\footnotetext[2]{sanjeev\_kumar\_19@yahoo.co.in}
\newcommand{\DS}{\c{S}tef\u{a}nescu }
\begin{abstract}
\c{S}tef\u{a}nescu  proved an elegant factorization result for polynomials over discrete valuation domains [CASC'2014, Lecture Notes in Computer Science, Ed. by V. Gerdt, W. Koepf, W. Mayr, and E. Vorozhtsov, Springer, Berlin, {Vol. \textbf{8660}}, pp. 460--471, 2014.] In this paper, a generalization of \c{S}tef\u{a}nescu's result is proved to cover a larger class of polynomials over discrete valuation domains. Such  results are useful in devising algorithms for polynomial factorization.
\end{abstract}
\section{Introduction}
Let $(R,v)$ be a discrete valuation domain. Let $f=a_0+a_1x+\cdots+a_nx^n \in {R}[x]$ be a nonconstant polynomial. Newton polygon $N_f$
of the polynomial $f$ is defined as the lower convex hull of the set $\{(i, v(a_i)) ~|~ a_i\neq 0\}$. The slopes of the underlying Newton polygon are the slopes of some line segments. Note that the slope of the line joining the points $(n, v(a_n))$ and
$(i, v(a_i))$ is $m_i(f)=(v(a_n)-v(a_i))/(n-i)$ for each $i=0, 1,\ldots,n-1$. The Newton index $e(f)$ of the polynomial $f$ is
 defined as
 \begin{eqnarray*}
e(f) = \max_{1\leq i\leq n}m_i(f).
 \end{eqnarray*}
It follows from the definition of $N_f$ and $e(f)$ that for nonconstant polynomials $f, g\in R[x]$, one has $e(fg) = \max(e(f), e(g))$. From the application point of view, Newton index has been used in devising algorithms for factoring polynomials \cite{S2}.
As a generalization of the classical result of Dumas \cite{Dumas}, \DS \cite{S1} proved a factorization result for polynomials over a discrete valuation domain using Newton index. Further, using the method of \cite{S1}, Kumar and Singh \cite{SKJS2022} extended the result of \DS to include a wider class of polynomials over discrete valuation domains. In \cite{S2}, \DS proved the following elegant factorization results.
\begin{thmx}\label{A}
Let $(R, v)$ be a discrete valuation domain. Let $f=a_0+a_1x+\cdots+a_nx^n \in {R}[x]$ be a nonconstant polynomial with $a_0a_n\neq 0$ and $n\geq 2$. Assume that there exists an index $s\in \{0, 1, 2, \ldots, n-1\}$ for which each of the following conditions is satisfied.
\begin{enumerate}[label=(\alph*)]
\item  $m_i(f)<m_s(f)$ for all $i\in \{0, 1, 2, \ldots, n-1\},~i\neq s$,
\item  $n(n-s)(m_s(f) -m_0(f))=1$,
\item $\gcd(v(a_s)-v(a_n), n-s)=1$.
\end{enumerate}
Then the polynomial $f$ is either irreducible in $R[x]$, or $f$ has a factor whose degree
is a multiple of $n-s$.
\end{thmx}
\begin{thmx}\label{B}
Let $(R, v)$ be a discrete valuation domain. Let $f=a_0+a_1x+\cdots+a_nx^n \in {R}[x]$ be a nonconstant polynomial with $a_0a_n\neq 0$ and $n\geq 2$. Assume that there exists an index $s\in \{0, 1, 2, \ldots, n-1\}$ for which each of the following conditions is satisfied.
\begin{enumerate}[label=(\alph*)]
\item  $m_i(f)<m_s(f)$ for all $i\in \{0, 1, 2,\ldots, n-1\},~i\neq s$,
\item  $u=n(n-s)(m_s(f) -m_0(f))\geq 2$,
\item  $\gcd(v(a_s)-v(a_n), n-s)=1$.
\end{enumerate}
Then either $f$ is irreducible in $R[x]$, or $f$ has a divisor whose degree is a multiple of $n-s$, or
 $f$ admits a factorization $f=f_1f_2$ such that $\alpha_2\deg(f_1)-\alpha_1 \deg(f_2)$ is a multiple of $n-s$  for some $\alpha_1, \alpha_2 \in \{1, \ldots, u-1\}$.
 \end{thmx}
This note ameliorates and extends the aforementioned factorization results on the lines of \cite{SKJS2022}. Our main results are the following:
\begin{theorem}\label{th1}
Let $(R, v)$ be a discrete valuation domain and let $f=a_0+a_1x+\cdots+a_nx^n \in {R}[x]$ be a nonconstant polynomial with $a_0a_n\neq 0$ and $n\geq 2.$ Assume that there exists an index $s\in \{0, 1, 2, \ldots, n-1\}$ such that the following conditions are satisfied:
\begin{enumerate}[label=(\alph*)]
\item\label{th1a}  $m_i(f)<m_s(f)$ for all $i\in \{0, 1, 2,\ldots, n-1\},~i\neq s$,
\item\label{th1b} $d=\gcd(v(a_s)-v(a_n), n-s)$ satisfies
\begin{eqnarray*}
d = \begin{cases} n(n-s)(m_s(f) -m_0(f)), &~\text{if}~s\neq 0;\\
1,  &~\text{if}~s=0.
\end{cases}
\end{eqnarray*}
\end{enumerate}
Then the polynomial $f$ is either irreducible in $R[x]$, or $f$ has a factor in $R[x]$ whose degree
is zero or a multiple of $(n-s)/d$.
\end{theorem}
\begin{theorem}\label{th2}
Let $(R, v)$ be a discrete valuation domain and let $f=a_0+a_1x+\cdots+a_nx^n \in {R}[x]$ be a nonconstant polynomial with $a_0a_n\neq 0$ and $n\geq 2$. Assume that there exists an index $s\in \{1, 2, \ldots, n-1\}$ such that each of the following conditions is satisfied.
\begin{enumerate}[label=(\alph*)]
\item\label{th2a}  $m_i(f)<m_s(f)$ for all $i\in \{0, 1, 2,\ldots, n-1\},~i\neq s$,
\item\label{th2b} $d=\gcd(v(a_s)-v(a_n), n-s)$ satisfies the following:
\begin{eqnarray*}
d=\begin{cases}
\text{a proper divisor of } u, \text{ where } u=n(n-s)(m_s(f) -m_0(f))\geq 2, & \text{ if } s\neq 0;\\
1, & \text{ if } s=0.
\end{cases}
\end{eqnarray*}
\end{enumerate}
Then either $f$ is irreducible in $R[x]$, or $f$ has a divisor whose degree is zero or a multiple of $(n-s)/d$, or
$f$ admits a factorization $f=f_1f_2$ such that $\alpha_2 \deg(f_1)-\alpha_1 \deg(f_2)$ is a multiple of $(n-s)/d$ for some $\alpha_1, \alpha_2\in \{1, \ldots, (u/d)-1\}$ with $\alpha_1+\alpha_2=u/d$.
\end{theorem}
We observe that Theorems \ref{th1} and \ref{th2} reduce to Theorems \ref{A} and \ref{B}, respectively for $d=1$ and $s\neq 0$.

Further, Theorem \ref{th1} reduces to the main result of \cite{SKJS2022} incase $v(a_n)=0$.

In view of Theorem \ref{th1}, if we take  $d=2$ and $0<s<n/2$ so that $(n-s)>n/2$, then either $f$ is irreducible, or $f$ has a factor whose degree is a multiple of $(n-s)/2$, say $m(n-s)/2$ for some positive integer $m$.  If possible, suppose that $m\geq 4$, then we have $n>m(n-s)/2\geq mn/4\geq n$, which is absurd, and so, we must have $m<4$. So, either $f$ is irreducible, or $f$ has a factor whose degree is equal to one of $(n-s)/2$ and $n-s$.
\begin{example}
For a prime $p$, let $v=v_p$ denotes the $p$-adic valuation  on $\mathbb{Q}$. For $n\geq 5$, consider the polynomial
\begin{eqnarray*}
X &=& a_0+p^{n-4}a_1x+p^{n(n-3)-1}(a_2 x^2+a_n x^{n})\in \mathbb{Z}[x],
\end{eqnarray*}
where $a_0,a_1,a_2,a_n\in\{1,2,\ldots,p-1\}$.
Here, we have
\begin{eqnarray*}
m_2(f)&=&\frac{v_p(p^{n(n-3)-1}a_n)-v_p(p^{n(n-3)-1}a_2)}{n-2}=0,\\
m_1(f)&=&\frac{v_p(p^{n(n-3)-1}a_n)-v_p(p^{n-4}a_1)}{n-1}=n-3,\\
m_0(f)&=&\frac{v_p(p^{n(n-3)-1}a_n)-v_p(a_0)}{n-0}=n-3-\frac{1}{n},
\end{eqnarray*}
which shows that $e(X)=m_1(f)$, and so, $s=1$. Further, we have  $n(n-1)(m_1(f)-m_0(f))=n-1$, and $\gcd(v_p(p^{n-4}a_1)-v_p(p^{n(n-3)-1}a_n),n-1)=\gcd((n-1)(n-3),n-1)=n-1$.
By Theorem \ref{th1}, the polynomial $X$ is either irreducible, or $X$ has a factor whose degree is a multiple of $n-1$.
\end{example}
\begin{example}
For a prime $p$, let $v_p$ be the $p$-adic valuation on $\mathbb{Q}$. For a positive integer $d\geq 2$,  we consider the polynomial
\begin{eqnarray*}
Y_d &=& p^{d+1}+p^{d-1}x^{d+1}+x^{d(d+1)}\in \mathbb{Z}[x].
\end{eqnarray*}
Here $n=d(d+1)$, $a_0=p^{d+1}$, $a_{d+1}=p^{d-1}$, $a_n=1$, and $a_j=0$ for all $j\not\in \{0,d+1,d(d+1)\}$. So, we have
\begin{eqnarray*}
m_{d+1}(f)&=&\frac{v_p(a_n)-v_p(a_{d+1})}{d(d+1)-(d+1)}=-\frac{1}{d+1},\\
m_0(f)&=&\frac{v_p(a_n)-v_p(a_0)}{d(d+1)-0}=-\frac{1}{d},
\end{eqnarray*}
which shows that $e(X)=m_{d+1}(f)$, and so, $s=d+1$. Further, $u=(d-1)(d+1)\geq 2$, since $d\geq 2$. Furthermore,  $\gcd(v_p(a_s)-v_p(a_n),n-s)=\gcd(d-1,d(d+1)-d-1)=d-1$, which divides $u$. Thus, by Theorem \ref{th2}, the polynomial $Y_d$ is irreducible, or has a factor whose degree is a multiple of $d+1$, or $f$ has a factorization $f=f_1f_2$ in $\mathbb{Z}[x]$ such that $\alpha_2\deg f_1-\alpha_1 \deg f_2$ is a multiple of $d+1$, for some $\alpha_1,\alpha_2 \in \{1,2,\ldots,d\}$ with $\alpha_1+\alpha_2=d+1$.
\end{example}
\section{Proof of Theorems \ref{th1} and \ref{th2}}
\begin{proof}[\bf Proof of Theorem \ref{th1}]
Our method of proof is similar to that of \cite{SKJS2022}. If $s=0$, then the Newton polygon of $f$ is a straight line segment joining the points $(0,v(a_0))$  and $(n,v(a_n))$, and so, by the classical result of Dumas \cite{Dumas} on factorization of polynomials via Newton polygons, it follows that $f$ is either irreducible, or $f$ has a factor of degree zero.

Now assume that $s>0$. Suppose that $f$ is not irreducible in $R[x]$ so that $f$ admits a factorization
$f = f_1f_2$ in $R[x]$ with $\min\{\deg f_1,f_2\}\geq 1$. For each $i=1,2$, let $n_i=\deg(f_i)$ so that
$n=n_1+n_2$, and
$f_i=\sum_{j=0}^{n_i}a_{ij}x^{j}$.
Consequently, we have $a_0=a_{10}a_{20}$ and $a_n=a_{1n_1}a_{2n_2}$ so that $v(a_0)=v(a_{10})+v(a_{20})$ and $v(a_{n})=v(a_{1n_1})+v(a_{2n_2})$. If we let
\begin{eqnarray*}
c_s=v(a_n)-v(a_s),~c_n=v(a_n)-v(a_0),~c_{i0}=v(a_{in_i})-v(a_{i0}),~i=1,2,
\end{eqnarray*}
then it follows that $c_n=c_{10}+c_{20}$. By the hypothesis \ref{th1a} and the identity $e(f_1f_2) = \max(e(f_1), e(f_2))$, we have
\begin{eqnarray*}
\dfrac{c_s}{n-s}=\dfrac{v(a_n)-v(a_s)}{n-s}=e(f)\geq e(f_1)\geq m_0(f_1)=\dfrac{c_{10}}{n_1}.
\end{eqnarray*}
We then have $\dfrac{c_s}{n-s}-\dfrac{c_{10}}{n_1}\geq 0$, and so
$\dfrac{c_s}{d}n_1-\dfrac{n-s}{d}c_{10}\geq 0$. Note that from the hypothesis \ref{th1b},
$c_s/d$ and $(n-s)/d$ are both integers.
Since $e(f)\geq e(f_2)$, we must have $c_s n_2-(n-s)c_{20}\geq 0$ and $\dfrac{c_s}{d}n_2-\dfrac{(n-s)}{d}c_{20}\geq 0$. By the hypothesis \ref{th1b}, we have the following:
\begin{eqnarray}\label{e1}
1=\dfrac{c_s}{d}n-\dfrac{(n-s)}{d}c_n=\Big(\dfrac{c_s}{d}n_1-\dfrac{(n-s)}{d}c_{10}\Big)+\Big(\dfrac{c_s}{d}n_2-\dfrac{(n-s)}{d}c_{20}\Big),
\end{eqnarray}
which shows that one of  the nonnegative integers
$\dfrac{c_s}{d}n_1-\dfrac{(n-s)}{d}c_{10}$ and $\dfrac{c_s}{d}n_2-\dfrac{(n-s)}{d}c_{20}$ is zero.

First assume that  $\dfrac{c_s}{d}n_1-\dfrac{(n-s)}{d}c_{10}= 0$. Then from \eqref{e1}, we have $\dfrac{c_s}{d}n_1=\dfrac{(n-s)}{d}c_{10}$. Since $\gcd({c_s}/{d},{(n-s)}/{d})=1$, the integer ${(n-s)}/{d}$ must divide $n_1$. Similarly, if we assume that $\dfrac{c_s}{d}n_2=\dfrac{(n-s)}{d}c_{20}$, then it follows from \eqref{e1} that ${(n-s)}/{d}$ must divide $n_2$. This completes the proof.
\end{proof}
\begin{proof}[\bf Proof of Theorem \ref{th2}] Assume that $f=f_1f_2$ for some nonconstant polynomials $f_1$ and $f_2$ in $R[x]$. For $s=0$, Theorem \ref{th2} reduces to Theorem \ref{th1}. So, we assume that $s>0$.  We use the notation same as in the proof of Theorem \ref{th1} so that we arrive at the following:
\begin{eqnarray*}
\dfrac{c_s}{d}n - \dfrac{(n-s)}{d}c_{n} = \dfrac{u}{d};~\dfrac{c_s}{d}n_i-\dfrac{(n-s)}{d}c_{i0}\geq 0,~i=1,2,
\end{eqnarray*}
where $n=n_1+n_2, c_n=c_{10}+c_{20}$, and
\begin{eqnarray}\label{e2}
\Big(\dfrac{c_s}{d}n_1-\dfrac{(n-s)}{d}c_{10}\Big)+\Big(\dfrac{c_s}{d}n_2-\dfrac{(n-s)}{d}c_{20}\Big)=\dfrac{u}{d}.
\end{eqnarray}
If $\dfrac{c_s}{d}n_i-\dfrac{(n-s)}{d}c_{i0}= 0$ for any $i\in \{1,2\}$, then as in Theorem \ref{th1}, we deduce that the degree of a divisor of
$f$ must be divisible by ${(n-s)}/{d}$.

If $\dfrac{c_s}{d}n_1-\dfrac{(n-s)}{d}c_{10}= 1$, then $\dfrac{c_s}{d}n_2-\dfrac{(n-s)}{d}c_{20}=\dfrac{u}{d}-1$, and so, from \eqref{e2}, we have
\begin{eqnarray*}
\dfrac{c_s}{d}\bigl(n_2-\bigl(\dfrac{u}{d}-1\bigr)n_1\bigr)=\dfrac{(n-s)}{d}\bigl(c_{20}-\bigl(\dfrac{u}{d}-1\bigr)c_{10}\bigr),
\end{eqnarray*}
which in view of the fact that $c_s/d$ and $(n-s)/d$ are coprime, shows that ${(n-s)}/{d}$ divides $n_2-(\dfrac{u}{d}-1)n_1$.
More generally, if we let
$\alpha_i=\dfrac{c_s}{d}n_i-\dfrac{(n-s)}{d}c_{i0}$, $i=1,2$ with
$\alpha_1+\alpha_2={u}/{d}$, then using \eqref{e2} one has the following:
\begin{eqnarray*}
\dfrac{c_s}{d}(\alpha_2 n_1-\alpha_1 n_2)=\dfrac{(n-s)}{d}(\alpha_2 c_{10}-\alpha_1 c_{20}).
\end{eqnarray*}
This in view of the fact that $\gcd({(n-s)}/{d}, {c_s}/{d})=1$ shows that ${(n-s)}/{d}$ must divide $\alpha_2 n_1-\alpha_1 n_2$.  This completes the proof.
\end{proof}
\subsection*{Acknowledgments.}
The present research is supported by Science and Engineering Research Board(SERB), a
statutory body of Department of Science and Technology (DST), Govt. of India through the project grant
no. MTR/2017/000575 awarded to the second author under MATRICS Scheme.
\subsection*{Disclosure} The authors declare to have no competing interest.

\end{document}